\documentclass[a4paper,11pt]{amsart}

\usepackage{graphicx}
\usepackage{mathptmx}
\usepackage{amsmath}
\usepackage{amssymb}
\usepackage{enumitem}
\usepackage{xcolor}

\newmuskip\pFqmuskip

\newcommand*\pFq[6][8]{%
  \begingroup 
  \pFqmuskip=#1mu\relax
  \mathcode`=\string"8000
  \begingroup\lccode`\~=`\,
  \lowercase{\endgroup\let~}\pFqcomma
  F^{#2}_{#3}{\left(\genfrac..{0pt}{}{#4}{#5}\bigg|#6\right)}%
  \endgroup
}
\newcommand{\pFqcomma}{\mskip\pFqmuskip}

\newtheorem{theorem}{Theorem}[section]

\begin{document}

\title[Probabilistic derangement numbers and polynomials]{Probabilistic derangement numbers and polynomials}

\author{Taekyun  Kim}
\address{Department of Mathematics, Kwangwoon University, Seoul 139-701, Republic of Korea}
\email{tkkim@kw.ac.kr}
\author{Dae San  Kim }
\address{Department of Mathematics, Sogang University, Seoul 121-742, Republic of Korea}
\email{dskim@sogang.ac.kr}


\maketitle

\begin{abstract}
Let $Y$ be a random variable such that the moment generating function of $Y$ exists in a neighborhood of the origin. The aim of this paper is to study probabilistic versions of the derangement polynomials, the derangement polynomials of type 2 and the $r$-derangement numbers, namely the probabilistic derangement polynomials associated with $Y$, the probabilistic derangement polynomials of type 2 associated with $Y$ and the probabilistic $r$-derangement numbers associated $Y$, respectively. We derive some properties, explicit expressions, certain identities and recurrence relations for those polynomials and numbers.  In addition, we consider the special case that $Y$ is the gamma random variable with parameters $\alpha,\beta > 0$.
\end{abstract}
 
\section{Introduction} 
In combinatorics, a derangement is a permutation that has no fixed points. The number of derangements of an $n$ element  set is called the $n$th derangement number. For $0\le r\le n$, the $n$th $r$-derangement number is the number of derangements of an $n+r$ element set under the restriction that the first $r$ elements are in disjoint cycles. \par
Let $Y$ be a random variable satisfying the moment condition (see \eqref{11}). The aim of this paper is to study the probabilistic derangement polynomials associated with $Y$, the probabilistic derangement polynomials of type 2 associated with $Y$ and the probabilistic $r$-derangement numbers associated $Y$, as probabilistic versions of the derangement polynomials, the derangement polynomials of type 2 and the $r$-derangement numbers, respectively. We derive some properties, explicit expressions, certain identities and recurrence relations for those polynomials and numbers.  In addition, we consider the special case that $Y$ is the gamma random variable with parameters $\alpha,\beta > 0$.  \par
The outline of this paper is as follows. In Section 1, we recall the derangement numbers, the $r$-derangement numbers, the derangement polynomials $D_{n}(x)$ and the derangement polynomials of type 2. We remind the reader of the Fubini polynomials $F_{n}(x)$, the Stirling numbers of the first kind, the unsigned Stirling numbers of the first kind and the Stirling numbers of the second kind ${n \brace k}$. Assume that $Y$ is a random variable such that the moment generating function of $Y$,\,\, $E[e^{tY}]=\sum_{n=0}^{\infty}\frac{t^{n}}{n!}E[Y^{n}], \quad (|t| <r)$, exists for some $r >0$. Let $(Y_{j})_{j\ge 1}$ be a sequence of mutually independent copies of the random variable $Y$, and let $S_{k}=Y_{1}+Y_{2}+\cdots+Y_{k},\,\, (k \ge 1)$,\,\, with \, $S_{0}=0$. Then we recall the probabilistic Stirling numbers of the second kind associated with $Y$, the probabilistic Bell polynomials associated with $Y$, $\phi_{n}^{Y}(x)$, and  the gamma random variable with parameters $\alpha,\,\beta >0$. Section 2 is the main results of this paper. Let $(Y_{j})_{j \ge1},\,\, S_{k},\,\, (k=0,1,\dots)$ be as in the above. Then we first define the probabilistic derangement polynomials associated $Y$, $D_{n}^{Y}(x)$. In Theorem 2.1, we find two explicit expressions for $D_{n}^{Y}(x)$. We derive a recurrence relation for $D_{n}^{Y}(x)$ in Theorem 2.2. In Theorem 2.3, we express $\phi_{n}^{Y}(1-x)$ as a finite sum involving $D_{k}(x)$. We define the probabilistic Euler numbers $E_{n}^{Y}$ in a natural manner. In Theorem 2.4, we express the convolution $\sum_{m=0}^{n}\binom{n}{m}\phi_{m}^{Y}(1-x)E_{n-m}^{Y}$, as an infinite sum involving $D_{m}(x)$. We define the probabilistic $r$-derangement numbers $D_{n}^{(r, Y)}$. We derive finite sum expressions for $D_{n}^{(r, Y)}$ in Theorems 2.5 and 2.6. In Theorem 2.7, the $n$th moment of $Y$ is expressed as a finite sum involving $D_{k+r}^{(r,Y)}$. We express $D_{n+r}^{(r,Y)}$ as a finite sum involving $D_{n-l}^{Y},\,\,(0 \le l \le n)$ in Theorem 2.8. We define the probabilistic derangement polynomials of type 2, $d_{n}^{Y}(x)$. In Theorem 2.9, we find an explicit expression and a recurrence relation for $d_{n}^{Y}(x)$. In Theorem 2.10, we derive a finite sum identity which is given by $\sum_{m=0}^{n}d_{m}^{Y}(x){n \brace m}=\sum_{m=0}^{n}\sum_{l=0}^{m}(-1)^{l}{m \brace l}\binom{n}{m}E[Y^{l}]F_{n-m}(x)$. We deduce an explicit expression for $d_{n}^{Y}(x)$ in Theorem 2.11. Finally, when $Y \sim \Gamma(1,1)$, we derive an expression of $D_{n}^{Y}(x)$ as a finite sum involving $D_{l}(1-x),\,\,(0 \le l \le n)$, in Theorem 2.12. For the rest of this section, we recall the facts that are needed throughout this paper. \par

\vspace{0.1in}

For $n\ge 0$, the derangement numbers $D_{n},\ (n\ge 0)$ are given by 
\begin{equation}
D_{n}=n!\sum_{k=0}^{n}\frac{(-1)^{k}}{k!}=\sum_{k=0}^{n}\binom{n}{k}(n-k)!(-1)^{k},\quad (\mathrm{see}\ [5,6,8,10,18,19,21,24]). \label{1}	
\end{equation}
From \eqref{1}, we easily derive the following equation 
\begin{equation}
\frac{1}{1-t}e^{-t}=\sum_{n=0}^{\infty}D_{n}\frac{t^{n}}{n!},\quad (\mathrm{see}\ [18,19,21,24,29]).\label{2}
\end{equation}
For $0\le r\le n$, the $r$-derangement numbers $D_{n}^{(r)}$ are given by 
\begin{equation}
\sum_{n=0}^{\infty}D_{n}^{(r)}\frac{t^{n}}{n!}=\frac{t^{r}}{(1-t)^{r+1}}e^{-t},\quad (\mathrm{see}\ [18,19,21,24]).\label{3}
\end{equation}
The derangement polynomials are defined by 
\begin{equation}
\frac{e^{xt}}{1-t}e^{-t}=\sum_{n=0}^{\infty}D_{n}	(x)\frac{t^{n}}{n!},\quad (\mathrm{see}\ [11,18,19,21,24,30]),\label{4}
\end{equation}
and the derangement polynomials of type 2 are given by 
\begin{equation}
\frac{1}{1-xt}e^{-t}=\sum_{n=0}^{\infty}d_{n}(x)\frac{t^{n}}{n!},\quad (\mathrm{see}\ [18]).\label{5}
\end{equation}
Note that 
\begin{equation*}
D_{n}(0)=D_{n}\quad\mathrm{and}\quad d_{n}(1)=D_{n},\ (n\ge 0). 
\end{equation*}
It is well known that the Fubini polynomials are defined by 
\begin{equation}
F_{n}(x)=\sum_{k=0}^{n}{n \brace k}k!x^{k},\quad (n\ge 0),\quad (\mathrm{see}\ [6,16,17,20,23,30]),\label{6}
\end{equation}
where ${n \brace k}$ are the Stirling numbers of the second kind given by 
\begin{equation}
\frac{1}{k!}\big(e^{t}-1\big)^{t}=\sum_{n=k}^{\infty}{n \brace k}\frac{t^{n}}{n!},\quad (k\ge 0),\quad (\mathrm{see}\ [1-32]). \label{7}
\end{equation}
From \eqref{6}, we have 
\begin{equation}
\frac{1}{1-x(e^{t}-1)}=\sum_{n=0}^{\infty}F_{n}(x)\frac{t^{n}}{n!},\quad (\mathrm{see}\ [6,17,23]).\label{8}
\end{equation}
The Stirling numbers of the first kind are given by 
\begin{equation}
\frac{1}{k!}\big(\log(1+t)\big)^{k}=\sum_{n=k}^{\infty}S_{1}(n,k)\frac{t^{n}}{n!},\quad (\mathrm{see}\ [6-30]).\label{9}	
\end{equation}
For $n\ge k \ge 0$, the unsigned Stirling numbers of the first kind are defined by 
\begin{equation}
{n \brack k}=(-1)^{n-k}S_{1}(n,k),\quad (\mathrm{see}\ [6-30]).\label{10}
\end{equation}
Let $Y$ be a random variable such that the moment generating function of $Y$ 
\begin{equation}
E\big[e^{tY}\big]=\sum_{n=0}^{\infty}E[Y^{n}]\frac{t^{n}}{n!},\quad (|t|<r)\quad \textrm{exists for some $r>0$.}\label{11}
\end{equation}
Assume that $(Y_{j})_{j\ge 1}$ is a sequence of mutually independent copies of $Y$, and that $S_{k}=Y_{1}+Y_{2}+\cdots+Y_{k},\ (k\ge 1),\ S_{0}=0$, (see [2,3,14,31]). The probabilistic Stirling numbers of the second kind associated with random variable $Y$ are defined by 
\begin{equation*}
{n \brace k}_{Y}=\frac{1}{k!}\sum_{j=0}^{k}\binom{k}{j}(-1)^{k-j}E\big[S_{j}^{n}\big],\quad (n\ge 0),\quad (\mathrm{see}\ [2,3,14,22,32]). 
\end{equation*}
Equivalently, they are given by
\begin{equation}
\frac{1}{k!}(E[e^{tY}]-1)^{k}=\sum_{n=k}^{\infty} { n \brace k}_{Y} \frac{t^n}{n!}. \label{12}
\end{equation}
When $Y=1,\ {n\brace k}_{Y}={n\brace k},\ (n\ge k\ge 0)$. The probabilistic Bell polynomials associated with random variable $Y$ are defined by 
\begin{equation}
\phi_{n}^{Y}(x)=\sum_{k=0}^{n}{n \brace k}_{Y}x^{k},\quad (n\ge 0),\quad (\mathrm{see}\ [3,14,32]).\label{13}
\end{equation}
From \eqref{13}, we have 
\begin{equation}
e^{x(E[e^{Yt}]-1)}=\sum_{n=0}^{\infty}\phi_{n}^{Y}(x)\frac{t^{n}}{n!}. \label{14}
\end{equation}
A continuous random variable $X$ whose density function is given by 
\begin{displaymath}
f(x)=\left\{\begin{array}{ccc}
\frac{\beta}{\Gamma(\alpha)} e^{-\beta x}(\beta x)^{\alpha-1}, & \textrm{if $x\ge 0$,} \\
0, & \textrm{if $x<0$,}
\end{array}\right. 
\end{displaymath}
for some $\alpha, \beta>0$ is said to be the gamma random variable with parameters $\alpha,\beta$, which is denoted by $X\sim\Gamma(\alpha,\beta),$ (see [28,31]). \par 
In particular, if $X\sim\Gamma(1,1)$, then we have 
\begin{align*}
E\Big[e^{(X-1+p)t}\Big]&=\int_{0}^{\infty}e^{(x-1+p)t}f(x)dx=\int_{0}^{\infty}e^{(x-1+p)t}e^{-x}dx \\
&=\frac{e^{-t}}{1-t}e^{pt}=\sum_{n=0}^{\infty}D_{n}(p)\frac{t^{n}}{n!},\quad (t<1). 
\end{align*}
Thus, we get 
\begin{displaymath}
E\big[(X-1+p)^{n}\big]=D_{n}(p),\quad (n\ge 0). 
\end{displaymath}

\section{Probabilistic derangement numbers and polynomials associated with random variables} 
Let $(Y_{k})_{k\ge 1}$ be a sequence of mutually independent copies of the random variable $Y$, and let 
\begin{displaymath}
	S_{0}=0, \,\,S_{k}=Y_{1}+Y_{2}+\cdots+Y_{k},\quad (k\ge 1). 
\end{displaymath}
Now, we define the {\it{probabilistic derangement polynomials associated with $Y$}} by 
\begin{equation}
\frac{e^{xt}}{1-t}E\big[e^{-tY}\big]=\sum_{n=0}^{\infty}D_{n}^{Y}(x)\frac{t^{n}}{n!}. \label{15}
\end{equation}
In particular, for $x=0,\ D_{n}^{Y}=D_{n}^{Y}(0)$ are called the {\it{probabilistic derangement numbers associated with $Y$}}. \par 
When $Y=1$, we have $D_{n}^{Y}(x)=D_{n}(x)$, $(n\ge 0)$. From \eqref{15}, we have 
\begin{align}
\sum_{n=0}^{\infty}D_{n}^{Y}(x)\frac{t^{n}}{n!}&=\frac{E\big[e^{-tY}\big]}{1-t}e^{xt}=\sum_{l=0}^{\infty}D_{l}^{Y}\frac{t^{l}}{l!} \sum_{m=0}^{\infty}x^{m}\frac{t^{m}}{m!} \label{16} \\
&=\sum_{n=0}^{\infty}\sum_{l=0}^{n}\binom{n}{l}D_{l}^{Y}x^{n-l}\frac{t^{n}}{n!}.\nonumber
\end{align}
On the other hand, by \eqref{15}, we also have 
\begin{align}
\sum_{n=0}^{\infty}D_{n}^{Y}(x)\frac{t^{n}}{n!}&=\frac{E\big[e^{-tY}\big]}{1-t}e^{xt}=\frac{1}{1-t}E\big[e^{t(x-Y)}] \label{17} \\
&=\sum_{l=0}^{\infty}t^{l} \sum_{m=0}^{\infty}E[(x-Y)^{m}]\frac{t^{m}}{m!} \nonumber\\
&=\sum_{n=0}^{\infty}n!\sum_{m=0}^{n}\frac{E\big[(x-Y)^{m}\big]}{m!}\frac{t^{n}}{n!}.\nonumber
\end{align}
Thus, by comparing the coefficients on both sides of \eqref{16} and \eqref{17}, we obtain
the following theorem.
\begin{theorem}
For $n\ge 0$, we have 
\begin{displaymath}
D_{n}^{Y}(x)= n!\sum_{m=0}^{n}\frac{E\big[(x-Y)^{m}\big]}{m!}=\sum_{l=0}^{n}\binom{n}{l}D_{l}^{Y}x^{n-l}. 
\end{displaymath}
\end{theorem}
From \eqref{15}, we note that 
\begin{align}
E\big[e^{t(x-Y)}\big]&=(1-t)\sum_{n=0}^{\infty}D_{n}^{Y}(x)\frac{t^{n}}{n!} \label{19} \\
&=\sum_{n=1}^{\infty}\big(D_{n}^{Y}(x)-nD_{n-1}^{Y}(x)\big)\frac{t^{n}}{n!}+1.\nonumber
\end{align}
On the other hand, by Taylor expansion, we also get 
\begin{align}
E\big[e^{t(x-Y)}\big]&=\sum_{n=0}^{\infty}E\big[(x-Y)^{n}\big]\frac{t^{n}}{n!} \label{20}\\
&=1+\sum_{n=1}^{\infty}E\big[(x-Y)^{n}\big]\frac{t^{n}}{n!}.\nonumber
\end{align}
Therefore, by \eqref{19} and \eqref{20}, we obtain the following theorem. 
\begin{theorem}
For $n\ge 1$, we have 
\begin{displaymath}
D_{n}^{Y}(x)-nD_{n-1}^{Y}(x)=E\big[(x-Y)^{n}\big].
\end{displaymath}
In particular, for $x=0$, we get 
\begin{displaymath}
D_{n}^{Y}-nD_{n-1}^{Y}=(-1)^{n}E[Y^{n}]. 
\end{displaymath}
\end{theorem}
Replacing $t$ by $1-E[e^{Yt}]$ in \eqref{4}, we get 
\begin{align}
e^{(1-x)(E[e^{Yt}]-1)}&=\sum_{k=0}^{\infty}D_{k}(x)\frac{1}{k!}\big(1-E[e^{Yt}]\big)^{k}E[e^{Yt}] \label{21} \\
&=\sum_{k=0}^{\infty}D_{k}(x)(-1)^{k}\sum_{l=k}^{\infty}{l \brace k}_{Y}\frac{t^{l}}{l!}\sum_{m=0}^{\infty}E[Y^{m}]\frac{t^{m}}{m!}\nonumber \\
&=\sum_{l=0}^{\infty}\sum_{k=0}^{l}D_{k}(x)(-1)^{l}{l \brace k}_{Y}\frac{t^{l}}{l!} \sum_{m=0}^{\infty}E[Y^{m}]\frac{t^{m}}{m!}\nonumber \\
&=\sum_{n=0}^{\infty}\sum_{l=0}^{n}\sum_{k=0}^{l}\binom{n}{l}{l \brace k}_{Y}D_{k}(x)(-1)^{l}E[Y^{n-l}]\frac{t^{n}}{n!}.\nonumber
\end{align}
Therefore, by \eqref{14} and \eqref{21}, we obtain the following theorem. 
\begin{theorem}
For $n\ge 0$, we have 
\begin{displaymath}
\phi_{n}^{Y}(1-x)= \sum_{l=0}^{n}\sum_{k=0}^{l}\binom{n}{l}{l \brace k}_{Y}D_{k}(x)(-1)^{l}E[Y^{n-l}].
\end{displaymath}
\end{theorem}
The {\it{probabilistic Euler numbers}} are given by 
\begin{equation}
\frac{2}{E[e^{tY}]+1}=\sum_{n=0}^{\infty}E_{n}^{Y}\frac{t^{n}}{n!}. \label{23}	
\end{equation}
When $Y=1$, $E_{n}^{Y}=E_{n},\ (n\ge 0)$, where $E_{n}$ are the ordinary Euler numbers given by 
\begin{displaymath}
\frac{2}{e^{t}+1}=\sum_{n=0}^{\infty}E_{n}\frac{t^{n}}{n!},\quad (\mathrm{see}\ [30 ]). 
\end{displaymath}
Replacing $t$ by $-E[e^{tY}]$ in \eqref{4}, we have 
\begin{align}
\frac{1}{1+E[e^{Yt}]}e^{(1-x)E[e^{Yt}]}&=\sum_{m=0}^{\infty}D_{m}(x)\frac{(-1)^{m}}{m!}\Big(E[e^{Yt}]\Big)^{m} \label{24} \\
&=\sum_{m=0}^{\infty}D_{m}(x)\frac{(-1)^{m}}{m!}E\Big[e^{(Y_{1}+\cdots+Y_{m})t}\Big] \nonumber \\
&=\sum_{n=0}^{\infty}\sum_{m=0}^{\infty}D_{m}(x)\frac{(-1)^{m}}{m!}E[S_{m}^{n}]\frac{t^{n}}{n!}.\nonumber	
\end{align}
On the other hand, by \eqref{23}, we get 
\begin{align}
\frac{1}{1+E[e^{Yt}]}e^{(1-x)E[e^{Yt}]}&=\frac{e^{1-x}}{2}\frac{2}{1+E[e^{Yt}]} e^{(1-x)(E[e^{Yt}]-1)}\label{25} \\
	&=\frac{e^{1-x}}{2}\sum_{l=0}^{\infty}E_{l}^{Y}\frac{t^{l}}{l!}\sum_{m=0}^{\infty}\phi_{m}^{Y}(1-x)\frac{t^{m}}{m!} \nonumber \\
	&=\frac{e^{1-x}}{2}\sum_{n=0}^{\infty}\sum_{m=0}^{n}\binom{n}{m}\phi_{m}^{Y}(1-x)E_{n-m}^{Y}\frac{t^{n}}{n!}. \nonumber
\end{align}
Therefore, by \eqref{24} and \eqref{25}, we obtain the following theorem. 
\begin{theorem}
For $n\ge 0$, we have 
\begin{displaymath}
\sum_{m=0}^{n}\binom{n}{m}\phi_{m}^{Y}(1-x)E_{n-m}^{Y}=2e^{x-1}\sum_{m=0}^{\infty}\frac{D_{m}(x)}{m!}(-1)^{m}E[S_{m}^{n}]. 
\end{displaymath}
\end{theorem}
For $0\le r\le n$, we consider the {\it{probabilistic $r$-derangement numbers associated with $Y$}}, which are defined by 
\begin{equation}
\frac{t^{r}}{(1-t)^{r+1}}E\big[e^{-Yt}\big]=\sum_{n=r}^{\infty}D_{n}^{(r,Y)}\frac{t^{n}}{n!}.\label{26}
\end{equation}
When $Y=1$, $D_{n}^{(r,Y)}=D_{n}^{(r)},\ (n\ge 0)$. From \eqref{26}, we have 
\begin{align}
\frac{t^{r}}{(1-t)^{r+1}}E\big[e^{-Yt}\big]&=\sum_{k=0}^{\infty}\binom{k+r}{r}t^{k+r}\sum_{m=0}^{\infty}(-1)^{m}E[Y^{m}]\frac{t^{m}}{m!}\label{27} \\
&=\sum_{k=r}^{\infty}\binom{k}{r}t^{k}\sum_{m=0}^{\infty}(-1)^{m}E[Y^{m}]\frac{t^{m}}{m!}\nonumber \\
&=\sum_{n=r}^{\infty}n!\sum_{k=r}^{n}\binom{k}{r}(-1)^{n-k}\frac{E[Y^{n-k}]}{(n-k)!}\frac{t^{n}}{n!}.\nonumber 	
\end{align}
Therefore, by \eqref{26} and \eqref{27}, we obtain the following theorem. 
\begin{theorem}
For $n\ge r\ge 0$, we have 
\begin{displaymath}
D_{n}^{(r,Y)}= n!\sum_{k=r}^{n}\binom{k}{r}(-1)^{n-k}\frac{E[Y^{n-k}]}{(n-k)!}.
\end{displaymath}
\end{theorem}
From \eqref{26}, we note that 
\begin{align}
\sum_{n=r}^{\infty}D_{n}^{(r,Y)}\frac{t^{n}}{n!}&=\bigg(\frac{t}{1-t}\bigg)^{r} \frac{1}{1-t}E[e^{-Yt}] \label{28}\\
&=\sum_{l=r}^{\infty}\binom{l-1}{r-1}t^{l}\sum_{m=0}^{\infty}D_{m}^{Y}\frac{t^{m}}{m!}\nonumber \\
&=\sum_{n=r}^{\infty}n!\sum_{l=r}^{n}\binom{l-1}{r-1}\frac{D_{n-l}^{Y}}{(n-l)!}\frac{t^{n}}{n!}.\nonumber
\end{align}
By comparing the coefficients on both sides of \eqref{28}, we obtain the following theorem. 
\begin{theorem}
For $n\ge r\ge 0$, we have 
\begin{displaymath}
D_{n}^{(r,Y)}= n!\sum_{l=r}^{n}\binom{l-1}{r-1}\frac{D_{n-l}^{Y}}{(n-l)!}.
\end{displaymath}
\end{theorem}
By \eqref{26}, we get 
\begin{align}
E\big[e^{-Yt}\big]&=(1-t)^{r+1}\sum_{k=r}^{\infty}D_{k}^{(r,Y)}\frac{t^{k-r}}{k!}\label{29}\\
&=(1-t)^{r+1}\sum_{k=0}^{\infty}D_{k+r}^{(r,Y)}\frac{t^{k}}{(k+r)!}\nonumber\\
&=\sum_{l=0}^{\infty}(-1)^{l}\binom{r+1}{l}t^{l}\sum_{k=0}^{\infty}D_{k+r}^{(r,Y)}\frac{t^{k}}{(k+r)!}\nonumber \\
&=\sum_{n=0}^{\infty}n!\sum_{k=0}^{n}\frac{D_{k+r}^{(r,Y)}}{(k+r)!}(-1)^{n-k}\binom{r+1}{n-k}\frac{t^{n}}{n!}.\nonumber	
\end{align}

Therefore, by \eqref{29}, we obtain the following theorem. 
\begin{theorem}
For $n\ge 0$, we have 
\begin{displaymath}
E\big[Y^{n}\big]=n!\sum_{k=0}^{n}\frac{D_{k+r}^{(r,Y)}}{(k+r)!}(-1)^{k}\binom{r+1}{n-k}.
\end{displaymath}
\end{theorem}
Now, we observe that 
\begin{align}
\frac{1}{(1-t)^{r+1}}E\big[e^{-Yt}\big]&=\frac{1}{t^{r}}\frac{t^{r}}{(1-t)^{r+1}}E\big[e^{-tY}\big]\label{31} \\
&=\sum_{n=0}^{\infty}\frac{D_{n+r}^{(r,Y)}}{(n+r)!}t^{n}=\sum_{n=0}^{\infty}\frac{D_{n+r}^{(r,Y)}}{\binom{n+r}{r}r!}\frac{t^{n}}{n!}. \nonumber
\end{align}
On the other hand, by binomial expansion, we get 
\begin{align}
\frac{1}{(1-t)^{r+1}}E\big[e^{-Yt}\big]&=\bigg(\frac{1}{1-t}\bigg)^{r}\frac{1}{1-t} E\big[e^{-Yt}\big]\label{32} \\
&=\sum_{l=0}^{\infty}\binom{r+l-1}{l}t^{l}\sum_{m=0}^{\infty}D_{m}^{Y}\frac{t^{m}}{m!}
\nonumber\\
&=\sum_{n=0}^{\infty}n!\sum_{l=0}^{n}\binom{r+l-1}{l}\frac{D_{n-l}^{Y}}{(n-l)!}\frac{t^{n}}{n!}.\nonumber	
\end{align}
Thus, by \eqref{31} and \eqref{32}, we get the following theorem. 
\begin{theorem}
For $n\ge 0$, we have 
\begin{displaymath}
D_{n+r}^{(r,Y)}=(n+r)!\sum_{l=0}^{n}\binom{r+l-1}{l}\frac{D_{n-l}^{Y}}{(n-l)!}. 
\end{displaymath}
\end{theorem}
Now, we define the {\it{probabilistic derangement polynomials of type 2 associated with $Y$}} as 
\begin{equation}
\frac{1}{1-xt}E\big[e^{-Yt}\big]=\sum_{n=0}^{\infty}d_{n}^{Y}(x)\frac{t^{n}}{n!}.\label{34}
\end{equation}
When $Y=1$, we have $d_{n}^{Y}(x)=d_{n}(x),\ (n\ge 0)$. In particular, for $x=1$, $d_{n}^{Y}(1)=D_{n}^{Y},\ (n\ge 0)$. \par 
From \eqref{34}, we have 
\begin{align}
\sum_{n=0}^{\infty}d_{n}^{Y}(x)\frac{t^{n}}{n!}&=\frac{1}{1-xt}E[e^{-Yt}] \label{35} \\
&=\sum_{m=0}^{\infty}x^{m}t^{m}\sum_{k=0}^{\infty}\frac{(-1)^{k}}{k!}E[Y^{k}]t^{k}\nonumber\\
&=\sum_{n=0}^{\infty}n!\sum_{k=0}^{n}\frac{(-1)^{k}}{k!}E\big[Y^{k}\big]x^{n-k}\frac{t^{n}}{n!},\nonumber
\end{align}
and 
\begin{align}
\sum_{n=0}^{\infty}(-1)^{n}E[Y^{n}]\frac{t^{n}}{n!}&=E\big[e^{-Yt}\big]=(1-xt)\sum_{n=0}^{\infty}d_{n}^{Y}(x)\frac{t^{n}}{n!} \label{36}\\
&=d_{0}^{Y}(x)+\sum_{n=1}^{\infty}\bigg(d_{n}^{Y}(x)-nxd_{n-1}^{Y}(x)\bigg)\frac{t^{n}}{n!}.\nonumber
\end{align}
Thus, by \eqref{36}, we get 
\begin{equation}
d_{0}^{Y}(x)=1,\quad d_{n}^{Y}(x)=nxd_{n-1}^{Y}(x)+(-1)^{n}E[Y^{n}],\quad (n\ge 1). \label{37}
\end{equation}
Therefore, by \eqref{35} and \eqref{37}, we obtain the following theorem. 
\begin{theorem}
For $n\ge 0$, we have 
\begin{displaymath}
d_{n}^{Y}(x)=n!\sum_{k=0}^{n}\frac{(-1)^{k}}{k!}E[Y^{k}]x^{n-k}. 
\end{displaymath}
Moreover, we have
\begin{displaymath}
d_{0}^{Y}(x)=1,\quad d_{n}^{Y}(x)=nxd_{n-1}^{Y}(x)+(-1)^{n}E[Y^{n}],\quad (n\ge 1). 
\end{displaymath}
\end{theorem}
Replacing $t$ by $e^{t}-1$ in \eqref{34}, we get 
\begin{align}
\frac{1}{1-x(e^{t}-1)}E\big[e^{-Y(e^{t}-1)}\big]&=\sum_{m=0}^{\infty}d_{m}^{Y}(x)\frac{1}{m!}(e^{t}-1)^{m} \label{38}\\
&=\sum_{n=0}^{\infty}\sum_{m=0}^{n}	d_{m}^{Y}(x){n \brace m}\frac{t^{n}}{n!}. \nonumber
\end{align}
On the other hand, by \eqref{8}, we obtain
\begin{align} 
\frac{1}{1-x(e^{t}-1)}E\big[e^{-Y(e^{t}-1)}\big]&=\sum_{k=0}^{\infty}F_{k}(x)\frac{t^{k}}{k!}\sum_{l=0}^{\infty}(-1)^{l}E[Y^{l}]\frac{(e^{t}-1)^{l}}{l!}\label{39} \\
&=\sum_{k=0}^{\infty}F_{k}(x)\frac{t^{k}}{k!}\sum_{l=0}^{\infty}(-1)^{l}E[Y^{l}]\sum_{m=l}^{\infty}{m \brace l}\frac{t^{m}}{m!} \nonumber \\
&=\sum_{k=0}^{\infty}F_{k}(x)\frac{t^{k}}{k!}\sum_{m=0}^{\infty}\sum_{l=0}^{m}(-1)^{l}E[Y^{l}]{m \brace l}\frac{t^{m}}{m!}\nonumber\\
&=\sum_{n=0}^{\infty}\sum_{m=0}^{n}\sum_{l=0}^{m}(-1)^{l}\binom{n}{m}{m \brace l}E[Y^{l}]F_{n-m}(x)\frac{t^{n}}{n!}.\nonumber
\end{align}
Therefore, by \eqref{38} and \eqref{39}, we obtain the following theorem. 
\begin{theorem}
For $n\ge 0$, we have 
\begin{displaymath}
\sum_{m=0}^{n}d_{m}^{Y}(x){n \brace m}=\sum_{m=0}^{n}\sum_{l=0}^{m}(-1)^{l}{m \brace l}\binom{n}{m}E[Y^{l}]F_{n-m}(x). 
\end{displaymath}
\end{theorem}
From \eqref{34}, we note that 
\begin{align}
\sum_{n=0}^{\infty}d_{n}^{Y}(x)\frac{t^{n}}{n!}&=e^{-\log(1-xt)}E\big[e^{-Yt}\big]=\sum_{l=0}^{\infty}\frac{(-1)^{l}}{l!}\big(\log(1-xt)\big)^{l}E\big[e^{-Yt}\big] \label{40}\\
&=\sum_{j=0}^{\infty}\sum_{l=0}^{j}{j \brack l}x^{j}\frac{t^{j}}{j!}\sum_{m=0}^{\infty}(-1)^{m}E[Y^{m}]\frac{t^{m}}{m!} \nonumber \\
&=\sum_{n=0}^{\infty}\sum_{j=0}^{n}\sum_{l=0}^{j}{j \brack l}\binom{n}{j}x^{j}(-1)^{n-j}E[Y^{n-j}]\frac{t^{n}}{n!}. \nonumber
\end{align}
Therefore, by comparing the coefficients on both sides of \eqref{40}, we obtain the following theorem. 
\begin{theorem}
For $n\ge 0$, we have 
\begin{displaymath}
d_{n}^{Y}(x)=\sum_{j=0}^{n}\sum_{l=0}^{j}{j \brack l}\binom{n}{j}x^{j}(-1)^{n-j}E[Y^{n-j}]. 
\end{displaymath}
\end{theorem}
Let $Y\sim \Gamma(1,1)$. Then, for $ t > -1$, we have 
\begin{align}
E\big[e^{(x-Y)t}\big]&=e^{xt}\int_{0}^{\infty}e^{-yt}e^{-y}dy=\frac{1}{1+t}e^{t}e^{(x-1)t} \label{41} \\
&=\sum_{l=0}^{\infty}(-1)^{l}D_{l}(1-x)\frac{t^{l}}{l!}.\nonumber	
\end{align}
From \eqref{41}, we have 
\begin{align}
\sum_{n=0}^{\infty}D_{n}^{Y}(x)\frac{t^{n}}{n!}&=\frac{e^{xt}}{1-t}E\big[e^{-Yt}\big]=\frac{1}{1-t}\frac{1}{1+t}e^{t}e^{(x-1)t}\label{42}\\
&=\frac{1}{1-t}\sum_{l=0}^{\infty}(-1)^{l}D_{l}(1-x)\frac{t^{l}}{l!}=\sum_{n=0}^{\infty}n!\sum_{l=0}^{n}\frac{(-1)^{l}}{l!}D_{l}(1-x)\frac{t^{n}}{n!}.\nonumber
\end{align}
Therefore, by comparing the coefficients on both sides of \eqref{42}, we obtain the following theorem. 
\begin{theorem}
Let $Y\sim\Gamma(1,1)$. For $n\ge 0$, we have 
\begin{displaymath}
D_{n}^{Y}(x)=n!\sum_{l=0}^{n}\frac{(-1)^{l}}{l!}D_{l}(1-x). 
\end{displaymath}
\end{theorem}
\section{Conclusion} 
By means of generating functions we studied probabilistic versions of the derangement polynomials, the $r$-derangement numbers and the derangement polynomials of type 2, namely the probabilistic derangement polynomials $D_{n}^{Y}(x)$ associated with $Y$, the probabilistic $r$-derangement numbers $D_{n}^{(r,Y)}$ associated $Y$ and the probabilistic derangement polynomials of type 2 $d_{n}^{Y}(x)$ associated with $Y$. Here $Y$ is a random variable such that the moment generating function of $Y$ exists in a neighborhood of the origin. In more detail, we derived an explicit expression for $D_{n}^{Y}(x)$ (see Theorem 2.1) and that in terms of derangement polynomials $D_{l}(x)$ for the special case of $Y \sim \Gamma(1,1)$ (see Theorem 2.12). We deduced a recurrence relation for $D_{n}^{Y}(x)$ (see Theorem 2.2). We expressed $\phi_{n}^{Y}(1-x)$ and the convolution of that with the probabilistic Euler numbers $E_{n}^{Y}$ in terms of $D_{k}(x)$ (see Theorems 2.3, 2.4). We found explicit expressions for $D_{n}^{(r,Y)}$ (see Theorems 2.5, 2.6) and two identities involving those numbers (see Theorems 2.7, 2.8). We obtained explicit expressions for $d_{n}^{Y}(x)$ (see Theorems 2.9, 2.11), a recurrence relation for $d_{n}^{Y}(x)$ (see Theorem 2.9) and an identity involving those polynomials (see Theorem 2.10). \par
As one of our future projects, we would like to continue to study probabilistic versions of many special polynomials and numbers and to find their applications to physics, science and engineering as well as to mathematics.

\end{document}